\documentclass[12pt]{article}%
\usepackage{graphicx}
\usepackage[intlimits]{amsmath}
\usepackage{latexsym}
\usepackage{amsfonts}
\usepackage{amssymb}%
\makeatletter
\newfont{\footsc}{cmcsc10 at 8truept}
\newfont{\footbf}{cmbx10 at 8truept}
\newfont{\footrm}{cmr10 at 10truept}
\newtheorem{theorem}{Theorem}

\newtheorem{lemma}[theorem]{Lemma}

\newenvironment{proof}[1][Proof]{\noindent{\textbf {#1}  }}  {\hfill$\Box$\bigskip}

\begin{document}

\title{An extension of Maclaurin's inequality}
\author{Vladimir Nikiforov\\Department of Mathematical Sciences, University of Memphis, \\Memphis TN 38152, USA, email: \textit{vnkifrv@memphis.edu}}
\maketitle

\begin{abstract}
Let $G$ be a graph of order $n$ and clique number $\omega.$ For every
$\mathbf{x}=\left(  x_{1},\ldots,x_{n}\right)  \in\mathbb{R}^{n}$ and $1\leq
s\leq\omega,$ set
\[
f_{s}\left(  G,\mathbf{x}\right)  =\sum\left\{  x_{i_{1}}\cdots x_{i_{s}%
}:\left\{  i_{1},\ldots,i_{s}\right\}  \text{ is an }s\text{-clique of
}G\right\}  ,
\]
and let $\rho_{s}\left(  G,\mathbf{x}\right)  =f_{s}\left(  G,\mathbf{x}%
\right)  \binom{\omega}{s}^{-1}.$ We show that if $\mathbf{x}\geq0,$ then
\[
\rho_{1}\left(  G,\mathbf{x}\right)  \geq\rho_{2}^{1/2}\left(  G,\mathbf{x}%
\right)  \geq\cdots\geq\rho_{\omega}^{1/\omega}\left(  G,\mathbf{x}\right)  .
\]
This extends the inequality of Maclaurin ($G=K_{n}$) and generalizes the
inequality of Motzkin and Straus. In addition, if $\mathbf{x}>0,$ for every
$1\leq s<\omega$ we determine when $\rho_{s}^{1/s}\left(  G,\mathbf{x}\right)
=\rho_{s+1}^{1/\left(  s+1\right)  }\left(  G,\mathbf{x}\right)  $.

Letting $k_{s}\left(  G\right)  $ be the number of $s$-cliques of $G,$ we show
that the above inequality is equivalent to the combinatorial inequality
\[
\frac{k_{1}\left(  G\right)  }{\binom{\omega}{1}}\geq\left(  \frac
{k_{2}\left(  G\right)  }{\binom{\omega}{2}}\right)  ^{1/2}\geq\cdots
\geq\left(  \frac{k_{\omega}\left(  G\right)  }{\binom{\omega}{\omega}%
}\right)  ^{1/\omega}.
\]

These results summarize previous work of Motzkin and Straus, Khadzhiivanov,
S\'{o}s and Straus, Fisher and Ryan, and Petingi and Rodriguez.\smallskip
\bigskip



\end{abstract}

\section{Introduction and main results}

Our graph-theoretic notation follows \cite{Bol98}; in particular, all graphs
are defined on the vertex set $\left\{  1,2,\ldots,n\right\}  =\left[
n\right]  $ and $G\left(  n\right)  $ stands for a graph with $n$ vertices. We
write $\omega\left(  G\right)  $ for the size of the maximal clique of $G$ and
$K_{s}\left(  G\right)  $ for the set of $s$-cliques of $G;$ we set
$k_{s}\left(  G\right)  =\left\vert K_{s}\left(  G\right)  \right\vert $.

For any graph $G=G\left(  n\right)  ,$ vector $\mathbf{x}=\left(  x_{1}%
,\ldots,x_{n}\right)  \in\mathbb{R}^{n},$ and $1\leq s\leq\omega=\omega\left(
G\right)  ,$ set
\[
f_{s}\left(  G,\mathbf{x}\right)  =\sum\left\{  x_{i_{1}}\cdots x_{i_{s}%
}:\left\{  i_{1},\ldots,i_{s}\right\}  \in K_{s}\left(  G\right)  \right\}
\]
and let $\rho_{s}\left(  G,\mathbf{x}\right)  =f_{s}\left(  G,\mathbf{x}%
\right)  \binom{\omega}{s}^{-1}.$ The inequality of Maclaurin (see, e.g.,
\cite{HLP52}, p. 52) reads as: if $G=K_{n}$ and $\mathbf{x}\geq0,$ then%
\begin{equation}
\rho_{1}\left(  G,\mathbf{x}\right)  \geq\rho_{2}^{1/2}\left(  G,\mathbf{x}%
\right)  \geq\cdots\geq\rho_{\omega}^{1/\omega}\left(  G,\mathbf{x}\right)  .
\label{mainin}%
\end{equation}
As it turns out, this inequality is valid for any graph $G$ and any
$\mathbf{x}\geq0.$ Moreover, letting $\mathbf{x}$ to be the vector of all
ones, we obtain%
\begin{equation}
\frac{k_{1}\left(  G\right)  }{\binom{\omega}{1}}\geq\left(  \frac
{k_{2}\left(  G\right)  }{\binom{\omega}{2}}\right)  ^{1/2}\geq\cdots
\geq\left(  \frac{k_{\omega}\left(  G\right)  }{\binom{\omega}{\omega}%
}\right)  ^{1/\omega}. \label{cliqin}%
\end{equation}
In particular, this inequality implies a concise form of Tur\'{a}n's theorem
\cite{Tur41}
\[
k_{2}\left(  G\right)  \leq\binom{\omega}{2}\left(  \frac{n}{\omega}\right)
^{2},
\]
and, more generally, of Zykov's theorem \cite{Zyk49}%
\[
k_{s}\left(  G\right)  \leq\binom{\omega}{s}\left(  \frac{n}{\omega}\right)
^{s}\text{ \ \ \ for every }2\leq s\leq\omega.
\]

To begin with, note that (\ref{mainin}) is essentially best possible. Indeed,
taking an $\omega$-clique $R$ in $G$ and letting $x_{i}=1/\omega$ if $i\in R,$
and $x_{i}=0$ if $i\notin R$, all inequalities in (\ref{mainin}) become equalities.

Note also that the inequality $\rho_{1}\left(  G,\mathbf{x}\right)  \geq
\rho_{2}^{1/2}\left(  G,\mathbf{x}\right)  $ has been proved by Motzkin and
Straus \cite{MoSt65}, so (\ref{mainin}) is an extension of their result.

In \cite{Kha77} Khadzhiivanov gave an analytical proof of inequality
(\ref{mainin}) and thus of (\ref{cliqin}), but his result remained unnoticed;
somewhat later S\'{o}s and Straus \cite{SoSt82} gave an independent analytical
proof of (\ref{cliqin}). Unfortunately, their result also remained generally
unknown, and so, in 1992, Fisher and Ryan \cite{FiRy92}, apparently unaware of
the previous work came up with a purely combinatorial proof of inequality
(\ref{cliqin}). Next, Petingi and Rodriguez \cite{PeRo00}, unaware of
\cite{Kha77} and \cite{SoSt82}, essentially rediscovered Khadzhiivanov's proof
of (\ref{mainin}), but without establishing the cases of equality. More
recently, Eckhoff \cite{Eck04}, apparently ignoring all of his predecessors,
found exactly $\max$ $k_{r}\left(  G\right)  $ for given $k_{2}\left(
G\right)  $ and $\omega\left(  G\right)  ,$ thus solving partially a problem
of Erd\H{o}s; his bound is as precise as one can get, yet its main term is
given by (\ref{cliqin}).

It should be noted, however, that the argument of Khadzhiivanov contains a gap
and his statement of the cases of equality in (\ref{mainin}) is incorrect.
Below we give a complete analytical proof of (\ref{mainin}) and determine the
cases of equality.

At first glance inequality (\ref{cliqin}) seems weaker than (\ref{mainin}),
yet in some sense they are equivalent since (\ref{cliqin}) implies in turn
(\ref{mainin}); in particular, Tur\'{a}n's theorem implies Motzkin-Straus's
result. Indeed, since $f_{s}\left(  G,\mathbf{x}\right)  $ is continuous in
$\mathbf{x},$ it suffices to deduce (\ref{mainin}) for all $\mathbf{x}$ with
positive rational coordinates. Moreover, since $f_{s}\left(  G,\mathbf{x}%
\right)  $ is a homogenous polynomial of degree $s,$ that is to say,
\begin{equation}
f_{s}\left(  G,a\mathbf{x}\right)  =a^{s}f_{s}\left(  G,\mathbf{x}\right)
\text{ for all }a\geq0,\text{ }\mathbf{x}\geq0, \label{eq2}%
\end{equation}
it suffices to deduce (\ref{mainin}) for all $\mathbf{x}$ with positive
integral entries. Let $x_{1},\ldots,x_{n}$ be positive integers; for every
$v\in V\left(  G\right)  ,$ replace $v$ by a set $U_{v}$ of size $x_{v}$ and
for every $uv\in E\left(  G\right)  ,$ replace $uv$ by a complete bipartite
graph with vertex classes $U_{u}$ and $U_{v}.$ Write $G_{\mathbf{x}}$ for the
resulting graph and note that $\omega\left(  G_{\mathbf{x}}\right)  =r$ and
$f_{s}\left(  G,\mathbf{x}\right)  =k_{s}\left(  G_{\mathbf{x}}\right)  .$
Hence, applying (\ref{cliqin}) to the graph $G_{\mathbf{x}},$ we see that
(\ref{mainin}) holds for $G$ and $\mathbf{x}=\left(  x_{1},\ldots
,x_{n}\right)  ,$ as claimed.

Thus, inequality (\ref{mainin}) is an analytical result that can be proved by
combinatorial means. The idea of this equivalence is not new and can be traced
back at least to Sidorenko \cite{Sid87}.

\bigskip

\section{\textbf{Proof of inequality (\ref{mainin})}}

In view of (\ref{eq2}), to prove (\ref{mainin}) for every graph $G=G\left(
n\right)  $ and every $s\in\left[  \omega\left(  G\right)  -1\right]  $, it
suffices to find $\max f_{s+1}\left(  G,\mathbf{x}\right)  ,$ subject to
$f_{s}\left(  G,\mathbf{x}\right)  =1.$ Let
\[
\mathcal{S}_{s}\left(  G\right)  =\{\mathbf{x}:\mathbf{x}\in\mathbb{R}%
^{n},\text{ }\mathbf{x}\geq0\ \text{and }f_{s}\left(  G,\mathbf{x}\right)
=1\}
\]
and note that the set $\mathcal{S}_{s}\left(  G\right)  $ is closed; for
$s\geq2$ it is unbounded and therefore, non-compact.

Our proof is based on two lemmas, the first of which establishes that
$f_{s+1}\left(  G,\mathbf{x}\right)  $ attains a maximum on $\mathcal{S}%
_{s}\left(  G\right)  :$ for $s\geq2$ this fact is not obvious.

\begin{lemma}
\label{le1} For every $G=G\left(  n\right)  $ and $1\leq s<\omega\left(
G\right)  ,$ the function $f_{s+1}\left(  G,\mathbf{x}\right)  $ attains a
maximum on $\mathcal{S}_{s}\left(  G\right)  .$
\end{lemma}

\begin{proof}
The lemma is obvious for $s=1$ since $\mathcal{S}_{1}\left(  G\right)  $ is
compact, so we shall assume $s\geq2.$ Our proof is by induction on $n.$ Let
$n=s+1,$ i.e., $G=K_{s+1}$. For every $\mathbf{x}\in\mathcal{S}_{s}\left(
G\right)  ,$ the AM-GM inequality implies that
\[
f_{s+1}\left(  G,\mathbf{x}\right)  =x_{1}x_{2}\cdots x_{s+1}\leq\left(
\frac{x_{1}\cdots x_{s}+\cdots+x_{2}x_{3}\cdots x_{s+1}}{s+1}\right)
^{\left(  s+1\right)  /s}=\left(  s+1\right)  ^{-\left(  s+1\right)  /s}.
\]
On the other hand, letting $\mathbf{y}=\left(  s+1\right)  ^{-1/s}\left(
1,\ldots,1\right)  \in\mathbb{R}^{s+1},$ we see that $f_{s}\left(
G,\mathbf{y}\right)  =1$ and $f_{s+1}\left(  G,\mathbf{y}\right)  =\left(
s+1\right)  ^{-\left(  s+1\right)  /s}.$ Hence, the assertion holds for
$n=s+1;$ assume that the assertion holds for any graph with fewer than $n$ vertices.

Suppose first that $G$ has a vertex $v$ that is not contained in any
$s$-clique of $G.$ We clearly have
\[
f_{s}\left(  G-v,\left(  x_{2},\ldots,x_{s+1}\right)  \right)  =f_{s}\left(
G,\mathbf{x}\right)  =1
\]
and $\ $%
\[
f_{s+1}\left(  G,\mathbf{x}\right)  =f_{s+1}\left(  G-v,\left(  x_{2}%
,\ldots,x_{s+1}\right)  \right)  .
\]
Since, by the induction hypothesis, the assertion holds for the graph $G-v,$
it holds for $G$ as well. So we may and shall assume that each vertex of $G$
is contained in an $s$-clique.

For all $\mathbf{x}\in\mathcal{S}_{s}\left(  G\right)  $ and all $\left\{
i_{1},\ldots,i_{s}\right\}  \in K_{s}\left(  G\right)  ,$ we have $x_{i_{1}%
}\cdots x_{i_{s}}\leq f_{s}\left(  G,\mathbf{x}\right)  =1.$ Thus, $x_{i_{1}%
}\cdots x_{i_{s}}\leq1$ for every $\left(  s+1\right)  $-clique $\left\{
i_{1},\ldots,i_{s+1}\right\}  ,$ and consequently, $f_{s+1}\left(
G,\mathbf{x}\right)  \leq\binom{n}{s+1}.$ Set
\[
M=\sup_{\mathbf{x}\in\mathcal{S}_{s}\left(  G\right)  }f_{s+1}(G,\mathbf{x})
\]
and, for every $i\geq1,$ select $\mathbf{x}^{\left(  i\right)  }=\left(
x_{1}^{\left(  i\right)  },\ldots,x_{n}^{\left(  i\right)  }\right)
\in\mathcal{S}_{s}\left(  G\right)  $ so that $\lim_{i\rightarrow\infty
}f_{s+1}(G,\mathbf{x}^{\left(  i\right)  })=M.$

To finish the proof, we shall find $\mathbf{y}\in\mathcal{S}_{s}\left(
G\right)  $ with $f_{s+1}(\mathbf{y})=M.$ If, for every $t\in\left[  n\right]
,$ the sequence $\left\{  x_{t}^{\left(  i\right)  }\right\}  _{i=1}^{\infty}$
is bounded, then $\left\{  \mathbf{x}^{\left(  i\right)  }\right\}
_{i=1}^{\infty}$ has an accumulation point $\mathbf{x}_{0}\in\mathcal{S}%
_{s}\left(  G\right)  ,$ and so $f_{s+1}\left(  G,\mathbf{x}_{0}\right)  =M,$
completing the proof. Assume now that $\left\{  x_{t}^{\left(  i\right)
}\right\}  _{i=1}^{\infty}$ is unbounded for some $t\in\left[  n\right]  $. By
assumption, $t\in R$ for some $R\in K_{s}\left(  G\right)  ;$ let say
$R=\left\{  1,\ldots,s-1,t\right\}  .$ Assume that there exists $c>0$ such
that $x_{v}^{\left(  i\right)  }>c$ for all $v\in\left[  s-1\right]  ,$
$i\geq1.$ Hence, for all $i\geq1,$
\[
M\geq f_{s}(G,\mathbf{x}^{\left(  i\right)  })\geq x_{1}^{\left(  i\right)
}\cdots x_{s-1}^{\left(  i\right)  }x_{t}^{\left(  i\right)  }>c^{s-1}%
x_{t}^{\left(  i\right)  },
\]
a contradiction, since $\left\{  x_{t}^{\left(  i\right)  }\right\}
_{i=1}^{\infty}$ is unbounded. Therefore, for some $v\in\left[  s-1\right]  ,$
the sequence $\left\{  x_{v}^{\left(  i\right)  }\right\}  _{i=1}^{\infty}$
contains arbitrarily small terms; let say $v=1$. Note that, for all $i\geq1,$%
\begin{align}
f_{s+1}\left(  G,\mathbf{x}^{\left(  i\right)  }\right)   &  \leq
x_{1}^{\left(  i\right)  }f_{s}\left(  G,\mathbf{x}^{\left(  i\right)
}\right)  +f_{s+1}\left(  G-v,\left(  x_{2}^{\left(  i\right)  },\ldots
,x_{n}^{\left(  i\right)  }\right)  \right) \nonumber\\
&  =x_{1}^{\left(  i\right)  }+f_{s+1}\left(  G-v,\left(  x_{2}^{\left(
i\right)  },\ldots,x_{n}^{\left(  i\right)  }\right)  \right)  \label{eq1}%
\end{align}
and
\[
f_{s}\left(  G-v,\left(  x_{2}^{\left(  i\right)  },\ldots,x_{n}^{\left(
i\right)  }\right)  \right)  \leq f_{s}\left(  G,\left(  x_{1}^{\left(
i\right)  },\ldots,x_{n}^{\left(  i\right)  }\right)  \right)  =1
\]
By the induction hypothesis, the function $f_{s+1}\left(  G-v,\mathbf{x}%
\right)  $ attains its maximum on $\mathcal{S}_{s}\left(  G-v\right)  ,$ let
say at $\mathbf{y}=\left(  y_{1},\ldots,y_{n-1}\right)  \in\mathcal{S}%
_{s}\left(  G-v\right)  ,$ and so
\[
f_{s+1}\left(  G-v,\left(  x_{2}^{\left(  i\right)  },\ldots,x_{n}^{\left(
i\right)  }\right)  \right)  \leq f_{s+1}\left(  G-v,\mathbf{y}\right)  .
\]
Hence, in view of (\ref{eq1}), we have
\[
f_{s+1}\left(  G-v,\mathbf{y}\right)  \leq M\leq x_{1}^{\left(  i\right)
}+f_{s+1}\left(  G-v,\mathbf{y}\right)  .
\]
Since $x_{1}^{\left(  i\right)  }$ can be arbitrarily small, it follows that
$f_{s+1}\left(  G-v,\mathbf{y}\right)  =M,$ and so
\[
f_{s+1}\left(  G,\left(  0,y_{1},\ldots,y_{n-1}\right)  \right)
=f_{s+1}\left(  G-v,\mathbf{y}\right)  =M,
\]
completing the proof.
\end{proof}

The proof of the next lemma essentially is due to Khadzhiivanov \cite{Kha77},
however, he overlooked the necessity of Lemma \ref{le1}.

\begin{lemma}
\label{le2} Assume that $G=G\left(  n\right)  $ is a noncomplete graph, $1\leq
s<\omega\left(  G\right)  ,$ every vertex of $G$ is contained in some
$s$-clique, and $f_{s+1}\left(  G,\mathbf{x}\right)  $ attains a maximum,
subject to $\mathbf{x}\in\mathcal{S}_{s}\left(  G\right)  $ at some
$\mathbf{y}>0$. If $u,v$ are nonadjacent vertices of $G,$ then there exists
$\mathbf{z}=\left(  z_{1},\ldots,z_{n}\right)  \in\mathcal{S}_{s}\left(
G\right)  $ such that $f_{s+1}\left(  G,\mathbf{z}\right)  =f_{s+1}\left(
G,\mathbf{y}\right)  $ and $z_{u}=0.$
\end{lemma}

\begin{proof}
By symmetry we shall assume that $u=1,$ $v=2$. For every $1\leq k\leq
\omega\left(  G\right)  ,$ $\xi,$ $\eta,$ and $\mathbf{x}=\left(  x_{1}%
,\ldots,x_{n}\right)  ,$ we have
\begin{equation}
f_{k}\left(  G,\left(  x_{1}+\xi,x_{2}+\eta,\ldots,x_{n}\right)  \right)
=\xi\frac{\partial f_{k}(G,\mathbf{x})}{\partial x_{1}}+\eta\frac{\partial
f_{k}(G,\mathbf{x})}{\partial x_{2}}+f_{k}(G,\mathbf{x}). \label{eq3}%
\end{equation}
Since $f_{s+1}(G,\mathbf{x})$ attains a maximum at $\mathbf{y}$, subject to
$f_{s}(G,\mathbf{x})=1,$ by Lagrange's method, there exists $\lambda$ such
that $\partial f_{s+1}(G,\mathbf{y})/\partial x_{i}=\lambda\partial
f_{s}(G,\mathbf{y})/\partial x_{i}$ for all $i\in\left[  n\right]  .$ Setting
\[
\xi=-y_{1},\text{ }\eta=y_{1}\frac{\partial f_{s}(G,\mathbf{y})/\partial
x_{1}}{\partial f_{s}(G,\mathbf{y})/\partial x_{2}},
\]
we see that
\[
\xi\frac{\partial f_{s}(G,\mathbf{y})}{\partial x_{1}}+\eta\frac{\partial
f_{s}(G,\mathbf{y})}{\partial x_{2}}=0
\]
and
\begin{equation}
\xi\frac{\partial f_{s+1}(G,\mathbf{y})}{\partial x_{1}}+\eta\frac{\partial
f_{s+1}(G,\mathbf{y})}{\partial x_{2}}=\lambda\left(  \xi\frac{\partial
f_{s}(G,\mathbf{y})}{\partial x_{1}}+\eta\frac{\partial f_{s}(G,\mathbf{y}%
)}{\partial x_{2}}\right)  =0. \label{eq4}%
\end{equation}
Hence, equality (\ref{eq3}) with $k=s,$ implies that
\begin{align*}
f_{s}\left(  G,\left(  0,y_{2}+\eta,y_{3},\ldots,y_{n}\right)  \right)   &
=f_{s}\left(  G,\left(  y_{1}+\xi,y_{2}+\eta,y_{3},\ldots,y_{n}\right)
\right)  =\\
&  =\xi\frac{\partial f_{s}(G,\mathbf{y})}{\partial x_{1}}+\eta\frac{\partial
f_{s}(G,\mathbf{y})}{\partial x_{2}}+f_{s}\left(  G,\mathbf{y}\right)  =1,
\end{align*}
and so, $\mathbf{z}=\left(  0,y_{2}+\eta,y_{3},\ldots,y_{n}\right)
\in\mathcal{S}_{s}\left(  G\right)  .$ On the other hand, equality (\ref{eq3})
with $k=s+1$ and (\ref{eq4}) imply that
\[
f_{s+1}\left(  G,\mathbf{z}\right)  =\xi\frac{\partial f_{s+1}(G,\mathbf{y}%
)}{\partial x_{1}}+\eta\frac{\partial f_{s+1}(G,\mathbf{y})}{\partial x_{2}%
}+f_{s+1}(G,\mathbf{y})=f_{s+1}(G,\mathbf{y}),
\]
completing the proof.
\end{proof}

To prove (\ref{mainin}) we first find $\mathbf{y}\in\mathcal{S}_{s}\left(
G\right)  $ such that $f_{s+1}\left(  G,\mathbf{y}\right)  \geq f_{s+1}\left(
G,\mathbf{x}\right)  $ for all $\mathbf{x}\in\mathcal{S}_{s}\left(  G\right)
.$ Set
\[
R=\left\{  v:v\in V\left(  G\right)  ,\text{ }y_{v}>0\text{ and }v\text{ is
contained in an }s\text{-clique}\right\}  ;
\]
without loss of generality we may assume that $G=G\left[  R\right]  $.
Applying Lemma \ref{le2} iteratively (i.e., using induction on $n)$, we see
that there exists $\mathbf{y}\in\mathcal{S}_{s}\left(  G\right)  $ such that
$f_{s+1}\left(  G,\mathbf{y}\right)  \geq f_{s+1}\left(  G,\mathbf{x}\right)
$ for all $\mathbf{x}\in\mathcal{S}_{s}\left(  G\right)  $ and the set
\[
R=\left\{  v:v\in V\left(  G\right)  ,\text{ }y_{v}>0\right\}
\]
induces a complete graph in $G;$ let $r=\left\vert R\right\vert \leq\omega.$
Maclaurin's inequality implies that
\[
\left(  f_{s+1}\left(  K_{r},\mathbf{y}\right)  \right)  ^{1/\left(
s+1\right)  }\leq\binom{r}{s+1}^{1/\left(  s+1\right)  }\binom{r}{s}%
^{-1/s}\leq\binom{\omega}{s+1}^{1/\left(  s+1\right)  }\binom{\omega}%
{s}^{-1/s},
\]
and so $\rho_{s+1}^{1/\left(  s+1\right)  }\left(  G,\mathbf{x}\right)
\leq\rho_{s}^{1/s}\left(  G,\mathbf{x}\right)  $ for every $1\leq
s<\omega\left(  G\right)  $ and $\mathbf{x}\geq0,$ completing the proof of
(\ref{mainin}).$\hfill\square$\bigskip

\subsection{\textbf{Cases of equality in (\ref{mainin})}}

Let $G=G\left(  n\right)  $ be a graph, $1\leq s<\omega=\omega\left(
G\right)  $, $\mathbf{x}=\left(  x_{1},\ldots,x_{n}\right)  >0,$ and $\rho
_{s}\left(  G,\mathbf{x}\right)  $ be defined as above; set%
\[
R_{s}=\left\{  v:v\in V\left(  G\right)  ,\text{ }v\text{ is contained in some
}s\text{-clique}\right\}  .
\]

\begin{theorem}
The equality $\rho_{s+1}^{1/\left(  s+1\right)  }\left(  G,\mathbf{x}\right)
=\rho_{s}^{1/s}\left(  G,\mathbf{x}\right)  $ holds if and only if $R_{s}$
induces a complete $\omega$-partite graph and if $V_{1},\ldots,V_{\omega}$ are
the vertex classes of $G\left[  R_{s}\right]  ,$ then $\sum_{v\in V_{i}}%
x_{v}=\sum_{v\in V_{j}}x_{v}$ for all $1\leq i<j\leq\omega.$
\end{theorem}

\begin{proof}
Assume $\rho_{s+1}^{1/\left(  s+1\right)  }\left(  G,\mathbf{x}\right)
=\rho_{s}^{1/s}\left(  G,\mathbf{x}\right)  ,$ set $\overline{R}%
=\overline{G\left[  R_{s}\right]  }$ - the complement of the graph induced by
$R_{s}$; let $G_{1},\ldots,G_{r}$ be the components of $\overline{R}.$
Clearly, $r\leq\omega;$ first we shall prove that $r=\omega$. Assume for
simplicity that $f_{s}\left(  G,\mathbf{x}\right)  =1;$ hence $f_{s+1}\left(
G,\mathbf{x}\right)  =\binom{\omega}{s+1}\binom{\omega}{s}^{-\left(
s+1\right)  /s}$ and $\mathbf{x}\in\mathcal{S}_{s}\left(  G\right)  .$
Applying Lemma \ref{le2}, preserving the value of $f_{s+1}\left(
G,\mathbf{x}\right)  ,$ find a vector $\mathbf{y}\in\mathcal{S}_{s}\left(
G\right)  $ with zero coordinates for all but one vertex from each component.
Then, by Maclaurin's inequality, we see that
\[
\binom{\omega}{s+1}\binom{\omega}{s}^{-\left(  s+1\right)  /s}=f_{s+1}\left(
G,\mathbf{y}\right)  \leq\binom{r}{s+1}\binom{r}{s}^{-\left(  s+1\right)
/s},
\]
and so $r=\omega.$ Clearly $G_{1},\ldots,G_{r}$ are complete subgraphs of
$\overline{R},$ since otherwise $K_{\omega+1}\subset G;$ hence, $G\left[
R_{s}\right]  $ is a complete $\omega$-partite graph. Setting $z_{i}%
=\sum_{v\in V_{i}}x_{v}$ for every $i\in\left[  \omega\right]  ,$ we see that
\begin{align*}
\left(  \binom{\omega}{s+1}^{-1}\sum_{1\leq i_{1}<\cdots<i_{s+1}\leq\omega
}z_{i_{1}}\ldots z_{i_{s+1}}\right)  ^{1/\left(  s+1\right)  }  &  =\rho
_{s+1}^{1/\left(  s+1\right)  }\left(  G,\mathbf{x}\right)  =\rho_{s}%
^{1/s}\left(  G,\mathbf{x}\right) \\
&  =\left(  \binom{\omega}{s}^{-1}\sum_{1\leq i_{1}<\cdots<i_{s}\leq\omega
}z_{i_{1}}\ldots z_{i_{s}}\right)  ^{1/s}.
\end{align*}
It is known (see, e.g., \cite{HLP52}, p. 52) that equality holds in
\[
\rho_{s+1}^{1/\left(  s+1\right)  }\left(  K_{\omega},\left(  z_{1}%
,\ldots,z_{\omega}\right)  \right)  =\rho_{s}^{1/s}\left(  K_{\omega},\left(
z_{1},\ldots,z_{\omega}\right)  \right)
\]
if and only if $z_{1}=\cdots=z_{\omega}$. Hence, the necessity of the
condition is proved. The sufficiency is immediate.
\end{proof}

\paragraph*{\textbf{Acknowledgment}}

Thanks are due to Jen\"{o} Lehel for valuable remarks.

\end{document}